\pdfoutput=1
\RequirePackage{ifpdf}
\ifpdf
\documentclass[pdftex]{sigma}
\else
\documentclass{sigma}
\fi

\begin{document}

\allowdisplaybreaks

\renewcommand{\PaperNumber}{099}

\FirstPageHeading

\ShortArticleName{On the Number of Real Roots of the Yablonskii--Vorob'ev Polynomials}

\ArticleName{On the Number of Real Roots\\ of the Yablonskii--Vorob'ev Polynomials}

\Author{Pieter ROFFELSEN}

\AuthorNameForHeading{P.~Rof\/felsen}

\Address{Radboud Universiteit Nijmegen, IMAPP, FNWI,\\
Heyendaalseweg 135, 6525 AJ Nijmegen, The Netherlands}
\Email{\href{mailto:roffelse@science.ru.nl}{roffelse@science.ru.nl}}

\ArticleDates{Received August 14, 2012, in f\/inal form December 07, 2012; Published online December 14, 2012}

\Abstract{We study the real roots of the Yablonskii--Vorob'ev polynomials, which are special polynomials used to represent rational solutions of the second Painlev\'e equation.
It has been conjectured that the number of real roots of the $n$th Yablonskii--Vorob'ev polynomial equals $\left[\frac{n+1}{2}\right]$.
We prove this conjecture using an interlacing property between the roots of the Yablonskii--Vorob'ev polynomials.
Furthermore we determine precisely the number of negative and~the number of positive real roots of the $n$th Yablonskii--Vorob'ev polynomial.}


\Classification{34M55}

\section{Introduction}

In this paper we study the real roots of the Yablonskii--Vorob'ev polynomials $Q_n$ ($n\in\mathbb{N}$).
Yablonskii and~Vorob'ev found these polynomials while studying the hierarchy of rational solutions of the second Painlev\'e equation.
The Yablonskii--Vorob'ev polynomials satisfy the def\/ining dif\/ferential-dif\/ference equation
\begin{gather*}
Q_{n+1}Q_{n-1}=zQ_n^2-4\big(Q_n Q_n''-(Q_n')^2\big),
\end{gather*}
with $Q_0=1$ and~$Q_1=z$.

The Yablonskii--Vorob'ev polynomials $Q_n$ are monic polynomials of degree $\frac{1}{2}n(n+1)$, with integer coef\/f\/icients.
The f\/irst few are given in Table~\ref{tableQn}.
Yablonskii~\cite{yablonskii} and~Vorob'ev~\cite{vorob} expressed the rational solutions of the second Painlev\'e equation,
\begin{gather*}P_{\rm II}(\alpha): \  w''(z)=2w(z)^3+zw(z)+\alpha,
\end{gather*}
with complex parameter $\alpha$, in terms of logarithmic derivatives of the Yablonskii--Vorob'ev polynomials, as summerized in the following theorem:

\begin{table}[t]
\centering
\caption{}\label{tableQn}
\vspace{1mm}
\begin{tabular}{r@{\,}c@{\,}l}
\hline
\multicolumn{3}{c}{Yablonskii--Vorob'ev polynomials\tsep{1pt}\bsep{1pt}}\\
\hline
$Q_2$ & $=$ & $4+z^3$\tsep{1pt}\\
$Q_3$ & $=$ & $-80+20z^3+z^6$\\
$Q_4$ & $=$ & $z\big(11200+60z^6+z^9\big)$\\
$Q_5$ & $=$ & $-6272000-3136000z^3+78400z^6+2800z^9+140z^{12}+z^{15}$\\
$Q_6$ & $=$ & $-38635520000+19317760000z^3+1448832000z^6-17248000z^9+
627200z^{12}$\\
& & ${}+18480z^{15}+280z^{18}+z^{21}$\\
$Q_7$ & $=$ & $z\big({-}3093932441600000-49723914240000z^6-828731904000z^9+13039488000
z^{12}$\\
& & ${} +62092800z^{15}+5174400z^{18}+75600z^{21}+504z^{24}+z^{27}\big)$\\
$Q_8$ & $=$ & $-991048439693312000000-743286329769984000000z^3$\\
& & ${} +37164316488499200000
z^6+1769729356595200000z^9+126696533483520000z^{12}$\\
& & ${} +407736096768000
z^{15}-6629855232000z^{18}+124309785600z^{21}+2018016000z^{24}$\\
& & ${} +32771200z^{27}+240240z^{30}+840z^{33}+z^{36}$\bsep{1pt}\\ \hline
\end{tabular}
\end{table}

\begin{theorem}
\label{thmYV}
$P_{\rm II}(\alpha)$ has a~rational solution iff $\alpha=n\in\mathbb{Z}$.
For $n\in\mathbb{Z}$ the rational solution is unique and~if $n\geq1$, then it is equal to
\begin{gather*}
w_n=\frac{Q_{n-1}'}{Q_{n-1}}-\frac{Q_n'}{Q_n}.
\end{gather*}
The other rational solutions are given by $w_0=0$ and~for $n\geq1$, $w_{-n}=-w_n$.
\end{theorem}

In~\cite{roffelsen} we proved the irrationality of the nonzero real roots of the Yablonskii--Vorob'ev polynomials, in this article we determine precisely the number of real roots of these polynomials.
Clarkson~\cite{clarksonoverview} conjectured that the number of real roots of $Q_n$ equals $\left[\frac{n+1}{2}\right]$, where $[x]$ denotes the integer part of $x$ for real numbers~$x$.
In Section~\ref{sectionreal} we prove this conjecture and~obtain the following theorem, where $Z_n$ is def\/ined as the set of real roots of~$Q_n$ for $n\in\mathbb{N}$.
\begin{theorem}
\label{thmnumberrealroots}
For every $n\in\mathbb{N}$, the number of real roots of $Q_n$ equals
\begin{gather}
\label{numberrealroots}
\left|Z_n\right|=\left[\frac{n+1}{2}\right].
\end{gather}
Furthermore for $n\geq2$,
\begin{gather}\label{minmaxroots}
\min(Z_{n-1})> \min(Z_{n+1}),\qquad
\max(Z_{n-1})< \max(Z_{n+1}).
\end{gather}
\end{theorem}

The argument is inductive and~an important ingredient is the fact that the real roots of $Q_{n-1}$ and~$Q_{n+1}$ interlace, which is proven by Clarkson~\cite{clarksonoverview}.

Kaneko and~Ochiai~\cite{kaneko} found a~direct formula for the lowest degree coef\/f\/icients of the Yablonskii--Vorob'ev polynomials $Q_n$ depending on~$n$.
In particular the sign of $Q_n(0)$ can be determined for $n\in\mathbb{N}$.
In Section~\ref{sectionpositivenegative} we use this to determine precisely the number of positive and~the number of negative real roots of $Q_n$, which yields to the following theorem.
\begin{theorem}
\label{thmpositivenegative}
Let $n\in\mathbb{N}$, then the number of negative real roots of $Q_n$ is equal to
\begin{gather*}
\left|Z_n\cap(-\infty,0)\right|=\left[\frac{n+1}{3}\right].
\end{gather*}
The number of positive real roots of $Q_n$ is equal to
\begin{gather*}
\left|Z_n\cap(0,\infty)\right|=
\begin{cases}\left[\dfrac{n}{6}\right]&\text{if $n$ is even,}\\[1ex]
\left[\dfrac{n+3}{6}\right] & \text{if $n$ is odd.}
\end{cases}
\end{gather*}
\end{theorem}
As a~consequence, for every $n\in\mathbb{N}$, we can calculate the number of positive real poles of the rational solution $w_n$ with residue $1$ and~with residue $-1$, and~the number of negative real poles of the rational solution $w_n$ with residue $1$ and~with residue $-1$.

\section{Number of real roots}
\label{sectionreal}

Let $P$ and~$Q$ be polynomials with no common real roots.
We say that the real roots of $P$ and~$Q$ interlace if and~only if in between any two real roots of $P$, $Q$ has a~real root and~in between any two real roots of~$Q$,~$P$ has a~real root.
Throughout this paper we use the convention $\mathbb{N}=\{0,1,2,\ldots\}$ and~def\/ine $\mathbb{N}^*:=\mathbb{N}\setminus\{0\}$.

\begin{theorem}
\label{simpleroots}
For every $n\in\mathbb{N}$, $Q_n$ has only simple roots.
Furthermore for $n\geq1$, $Q_{n-1}$ and~$Q_{n+1}$ have no common roots and~$Q_{n-1}$ and~$Q_n$ have no common roots.
\end{theorem}

\begin{proof}
See Fukutani, Okamoto and~Umemura~\cite{fukutani}.
\end{proof}

\begin{theorem}
\label{thminterlace}
For every $n\geq1$, the real roots of $Q_{n-1}$ and~$Q_{n+1}$ interlace.
\end{theorem}
\begin{proof}
See Clarkson~\cite{clarksonoverview}.
\end{proof}

Let $f,g:\mathbb{R}\rightarrow\mathbb{R}$ be continuous functions and~$x\in\mathbb{R}$.
We say that $f$ crosses $g$ positively at $x$ if and~only if $f(x)=g(x)$ and~there is a~$\delta>0$ such that $f(y)<g(y)$ for $x-\delta<y<x$ and~$f(y)>g(y)$ for $x<y<x+\delta$.
We say that $f$ crosses $g$ negatively at $x$ if and~only if $f(x)=g(x)$ and~there is a~$\delta>0$ such that $f(y)>g(y)$ for $x-\delta<y<x$ and~$f(y)<g(y)$ for $x<y<x+\delta$.
So $f$ crosses $g$ negatively at $x$ if and~only if $g$ crosses $f$ positively at $x$.

Let $m\in\mathbb{N}$ and~suppose that $f$ is $m$ times dif\/ferentiable, then we denote the $m$th derivative of $f$ by $f^{(m)}$ with convention $f^{(0)}=f$.

\begin{proposition}
\label{propanalytic}
Let $f,g:\mathbb{R}\rightarrow\mathbb{R}$ be analytic functions and~$x\in\mathbb{R}$.
Then $f$ crosses~$g$ positively at~$x$ if and~only if there is an $m\geq1$ such that $f^{(i)}(x)=g^{(i)}(x)$ for $0\leq i<m$ and~$f^{(m)}(x)>g^{(m)}(x)$.

Similarly $f$ crosses $g$ negatively at $x$ if and~only if there is a~$m\geq1$ such that $f^{(i)}(x)=g^{(i)}(x)$ for $0\leq i<m$ and~$f^{(m)}(x)<g^{(m)}(x)$.
\end{proposition}

\begin{proof}
This is proven easily using Taylor's theorem.
\end{proof}
\begin{lemma}
For every $n\in\mathbb{N}^*$ we have
\begin{subequations}
\begin{gather}
Q_{n+1}'Q_{n-1}-Q_{n+1}Q_{n-1}' =(2n+1)Q_n^2,\label{eqfukutani1}\\
Q_{n+1}''Q_{n-1}-Q_{n+1}Q_{n-1}'' =2(2n+1)Q_nQ_n',\label{eqfukutani2}\\
Q_{n+1}'''Q_{n-1}-Q_{n+1}Q_{n-1}''' =2(2n+1)\left(Q_n'\right)^2+(2n+1)Q_nQ_n''.\label{eqfukutani3}
\end{gather}
\end{subequations}
\end{lemma}
\begin{proof}
See Fukutani, Okamoto and~Umemura~\cite{fukutani}.
\end{proof}
The following proposition contains some well-known properties of the Yablonskii--Vorob'ev polynomials, see for instance Clarkson and~Mansf\/ield~\cite{clarksonmansfield}.
\begin{proposition}
\label{proplimitbehaviour}
For every $n\in\mathbb{N}$, $Q_n$ is a~monic polynomial of degree $\tfrac{1}{2}n(n+1)$ with integer coefficients.
As a consequence, for $n\geq1$,
\begin{gather*}
\lim_{x\rightarrow \infty} Q_n(x)=\infty,\qquad
\lim_{x\rightarrow -\infty} Q_n(x)=
\begin{cases} -\infty & \text{if $n\equiv1,2\pmod{4}$,}\\
\infty & \text{if $n\equiv0,3\pmod{4}$.}
\end{cases}
\end{gather*}
\end{proposition}

By Proposition~\ref{proplimitbehaviour}, $Q_n$ has real coef\/f\/icients and~hence we can consider $Q_n$ as a~real-valued function def\/ined on the real line, that is, we consider
\begin{gather*}
Q_n: \ \mathbb{R}\rightarrow\mathbb{R}.
\end{gather*}
\begin{proposition}
\label{crossing}
Let $n\in\mathbb{N}^*$, if $x\in\mathbb{R}$ is such that $Q_{n+1}$ crosses $Q_{n-1}$ positively at $x$, then
\begin{gather*}
Q_{n+1}(x)=Q_{n-1}(x)>0.
\end{gather*}
Similarly if $x\in\mathbb{R}$ is such that $Q_{n+1}$ crosses $Q_{n-1}$ negatively at $x$, then
\begin{gather*}
Q_{n+1}(x)=Q_{n-1}(x)<0.
\end{gather*}
\end{proposition}
\begin{proof}
Let $n\in\mathbb{N}^*$.
Suppose $x\in\mathbb{R}$ is such that $Q_{n+1}$ crosses $Q_{n-1}$ positively at $x$.
If
\begin{gather*}
Q_{n+1}(x)=Q_{n-1}(x)=0,
\end{gather*}
then $Q_{n+1}$ and~$Q_{n-1}$ have a~common root, which contradicts Theorem~\ref{simpleroots}.

Let us assume
\begin{gather}
\label{negativecrossing}
Q_{n+1}(x)=Q_{n-1}(x)<0.
\end{gather}
Then by Proposition~\ref{propanalytic},
\begin{gather}
\label{derivativenegative}
Q_{n+1}'(x)-Q_{n-1}'(x)\geq0.
\end{gather}
Therefore, by equation~\eqref{eqfukutani1},
\begin{gather*}
0\leq (2n+1)Q_n(x)^2=Q_{n+1}'(x)Q_{n-1}(x)-Q_{n+1}(x)Q_{n-1}'(x)\\
\phantom{0\leq (2n+1)Q_n(x)^2}
=Q_{n+1}(x)\left(Q_{n+1}'(x)-Q_{n-1}'(x)\right)\leq 0,
\end{gather*}
where in the last inequality we used equation~\eqref{negativecrossing} and~equation~\eqref{derivativenegative}.

We conclude
\begin{gather*}
(2n+1)Q_n(x)^2=Q_{n+1}(x)\left(Q_{n+1}'(x)-Q_{n-1}'(x)\right)=0,
\end{gather*}
so $Q_n(x)=0$ and~$Q_{n+1}'(x)=Q_{n-1}'(x)$.
Therefore by equation~\eqref{eqfukutani2},
\begin{gather*}
Q_{n+1}(x)\left(Q_{n+1}''(x)-Q_{n-1}''(x)\right)=Q_{n+1}''(x)Q_{n-1}(x)-Q_{n+1}(x)Q_{n-1}''(x)\\
\phantom{Q_{n+1}(x)\left(Q_{n+1}''(x)-Q_{n-1}''(x)\right)}
=2(2n+1)Q_n(x)Q_n'(x)=0.
\end{gather*}
We conclude $Q_{n+1}''(x)=Q_{n-1}''(x)$.
Since $Q_n(x)=0$ and, by Theorem~\ref{simpleroots}, $Q_n$ has only simple roots, we have $Q_n'(x)\neq0$.
Therefore by~\eqref{eqfukutani3},
\begin{gather*}
Q_{n+1}(x)\left(Q_{n+1}'''(x)-Q_{n-1}'''(x)\right)=Q_{n+1}'''(x)Q_{n-1}(x)-Q_{n+1}(x)Q_{n-1}'''(x)\\
\phantom{Q_{n+1}(x)\left(Q_{n+1}'''(x)-Q_{n-1}'''(x)\right)}=2(2n+1)\left(Q_n'(x)\right)^2+(2n+1)Q_n(x)Q_n''(x)\\
\phantom{Q_{n+1}(x)\left(Q_{n+1}'''(x)-Q_{n-1}'''(x)\right)}=2(2n+1)\left(Q_n'(x)\right)^2>0.
\end{gather*}
Since $Q_{n+1}(x)<0$ we conclude $Q_{n+1}'''(x)<Q_{n-1}'''(x)$.
So $Q_{n+1}'(x)=Q_{n-1}'(x)$, $Q_{n+1}''(x)=Q_{n-1}''(x)$ but $Q_{n+1}'''(x)<Q_{n-1}'''(x)$.
Therefore by Proposition~\ref{propanalytic}, $Q_{n+1}$ does not cross $Q_{n-1}$ positively at $x$ and~we have obtained a~contradiction.
We conclude that
\begin{gather*}
Q_{n+1}(x)=Q_{n-1}(x)>0.
\end{gather*}
The second part of the proposition is proven similar.
\end{proof}
We prove theorem~\ref{thmnumberrealroots}, using Theorem~\ref{thminterlace} and~Proposition \ref{crossing}.

\begin{proof}[Proof of Theorem~\ref{thmnumberrealroots}.]
Observe that~\eqref{numberrealroots} is correct for $n=0,1,2,3,4$.
Furthermore it is easy to see that \eqref{minmaxroots} is true for $n=1,2,3$.
We proceed by induction, suppose $n\geq4$ and
\begin{gather*}
\left|Z_{n-1}\right|=\left[\frac{n}{2}\right].
\end{gather*}
Then $Q_{n-1}$ has at least $2$ real roots.
By Theorem~\ref{thminterlace} the real roots of $Q_{n-1}$ and~$Q_{n+1}$ interlace, hence $Q_{n+1}$ has a~real root.
Let us def\/ine
\begin{gather*}
z:=\min(Z_{n+1}), \qquad z_1:=\min(Z_{n-1}), \qquad z_2:=\min(Z_{n-1}\setminus\left\{z_1\right\}),
\end{gather*}
so $z$ is the smallest real root of $Q_{n+1}$ and~$z_1$ and~$z_2$ are the smallest and~second smallest real root of $Q_{n-1}$ respectively.

By Theorem~\ref{thminterlace} the real roots of $Q_{n-1}$ and~$Q_{n+1}$ interlace, hence either $z<z_1$ or $z_1<z<z_2$.
We prove that $z_1<z<z_2$ can not be the case.
Suppose $z_1<z<z_2$ and~suppose $n\equiv0,1\pmod{4}$, then by Proposition~\ref{proplimitbehaviour},
\begin{gather*}
\lim_{x\rightarrow-\infty}Q_{n-1}(x)=\infty.
\end{gather*}
Hence $Q_{n-1}(x)>0$ for $x<z_1$.
Since $Q_{n-1}(z_1)=0$, this implies $Q_{n-1}'(z_1)\leq0$.
By Theorem~\ref{simpleroots}, $Q_{n-1}$ has only simple roots, hence $Q_{n-1}'(z_1)\neq0$, so $Q_{n-1}'(z_1)<0$.
Therefore by Proposition \ref{propanalytic}, $Q_{n-1}$ crosses $0$ negatively at $z_1$.
Hence $Q_{n-1}(x)<0$ for $z_1<x<z_2$, in particular
\begin{gather}
\label{negativez}
Q_{n-1}(z)<0.
\end{gather}
Since $n\equiv0,1\pmod{4}$, we have by Proposition~\ref{proplimitbehaviour},
\begin{gather*}
\lim_{x\rightarrow-\infty}Q_{n+1}(x)=-\infty.
\end{gather*}
Therefore $Q_{n+1}(x)<0$ for $x<z$, in particular
\begin{gather*}
Q_{n+1}(z_1)<0.
\end{gather*}
Def\/ine the polynomial $P:=Q_{n+1}-Q_{n-1}$, then
\begin{gather*}
P(z_1)=Q_{n+1}(z_1)-Q_{n-1}(z_1)=Q_{n+1}(z_1)-0<0,
\end{gather*}
and by equation~\eqref{negativez},
\begin{gather*}
P(z)=Q_{n+1}(z)-Q_{n-1}(z)=0-Q_{n-1}(z)>0.
\end{gather*}
So $P$ is a~polynomial with $P(z_1)<0$, $P(z)>0$ and~$z_1<z$.
Hence there is a~$z_1<x<z$ such that $P$ crosses $0$ positively at $x$, for instance
\begin{gather*}
x:=\inf\left\{t\in(z_1,z)\mid P(t)>0\right\},
\end{gather*}
has the desired properties.

Since $P=Q_{n+1}-Q_{n-1}$ crosses $0$ positively at $x$, $Q_{n+1}$ crosses $Q_{n-1}$ positively at $x$.
But $z_1<x<z_2$, hence
\begin{gather}
\label{minimum}
Q_{n+1}(x)=Q_{n-1}(x)<0.
\end{gather}
This contradicts Proposition~\ref{crossing}.

If $n\equiv2,3\pmod{4}$, then by a~similar argument, there is a $z_1<x<z$ such that $Q_{n+1}$ cros\-ses~$Q_{n-1}$ negatively at $x$ with
\begin{gather*}
Q_{n+1}(x)=Q_{n-1}(x)>0,
\end{gather*}
which again contradicts Proposition~\ref{crossing}.

We conclude that $z_1<z<z_2$ can not be the case and~hence $z<z_1$, that is,
\begin{gather*}
\min(Z_{n-1})>\min(Z_{n+1}).
\end{gather*}
Let us def\/ine
\begin{gather*}
w:=\max(Z_{n+1}), \qquad w_1:=\max(Z_{n-1}), \qquad w_2:=\max(Z_{n-1}\setminus\left\{w_1\right\}),
\end{gather*}
so $w$ is the largest real root of $Q_{n+1}$ and~$w_1$ and~$w_2$ are the largest and~second largest real root of~$Q_{n-1}$ respectively.

Suppose $w_1>w$, then by a~similar argument as the above, there is a $w_2<x<w$ such that~$Q_{n+1}$ crosses $Q_{n-1}$ positively at $x$ with
\begin{gather*}
Q_{n+1}(x)=Q_{n-1}(x)<0.
\end{gather*}
This is in contradiction with Proposition~\ref{crossing}, so $w_1<w$, that is
\begin{gather}
\label{maximum}
\max(Z_{n-1})<\max(Z_{n+1}).
\end{gather}
Let $z_1<z_2<\cdots<z_k$ be the real roots of $Q_{n-1}$ with $k=\left[\frac{n}{2}\right]$ and~$z_1'<z_2'<\cdots <z_m'$ be the real roots of $Q_{n+1}$.
Then by equations~\eqref{minimum} and~\eqref{maximum}, $z_1'<z_1$, $z_k<z_m'$ and~since by Theorem~\ref{thminterlace} the real roots of $Q_{n-1}$ and~$Q_{n+1}$ interlace, we have
\begin{gather*}
z_1'<z_1<z_2'<z_2<z_3'<z_3<\cdots<z_{k-1}'<z_{k-1}<z_{k}'<z_k<z_{k+1}'=z_m'.
\end{gather*}
Hence $m=k+1$, that is,
\begin{gather*}
\left|Z_{n+1}\right|=m=k+1=\left[\frac{n}{2}\right]+1=\left[\frac{n+2}{2}\right].
\end{gather*}
The theorem follows by induction.
\end{proof}

\section{Number of positive and~negative real roots}
\label{sectionpositivenegative}
For a~polynomial $P$ we denote the set of real roots of $P$ by $Z_P$.
\begin{lemma}
\label{lemnumbersroots}
Let $P$ and~$Q$ be polynomials with real coefficients, both a~positive leading coefficient and~only simple roots.
Assume that the real roots of $P$ and~$Q$ interlace.
Furthermore suppose both $P$ and~$Q$ have a~real root and
\begin{gather*}
\min(Z_P)>\min(Z_Q),\qquad
\max(Z_P)<\max(Z_Q).
\end{gather*}
Then we have the following relations between the number of negative and~positive real roots of~$P$ and~$Q$,
\begin{gather*}
\left|Z_Q\cap (-\infty,0)\right|=\left|Z_P\cap (-\infty,0)\right|+
\begin{cases} 1 & \text{if $P(0)=0$,}\\
0 & \text{if $Q(0)=0$,}\\
1 & \text{if $P(0)>0$ and~$Q(0)>0$,}\\
0 & \text{if $P(0)>0$ and~$Q(0)<0$,}\\
0 & \text{if $P(0)<0$ and~$Q(0)>0$,}\\
1 & \text{if $P(0)<0$ and~$Q(0)<0$,}
\end{cases}\\
\left|Z_Q\cap (0,\infty)\right|=\left|Z_P\cap (0,\infty)\right|+
\begin{cases} 1 & \text{if $P(0)=0$,}\\
0 & \text{if $Q(0)=0$,}\\
0 & \text{if $P(0)>0$ and~$Q(0)>0$,}\\
1 & \text{if $P(0)>0$ and~$Q(0)<0$,}\\
1 & \text{if $P(0)<0$ and~$Q(0)>0$,}\\
0 & \text{if $P(0)<0$ and~$Q(0)<0$.}
\end{cases}
\end{gather*}
\end{lemma}

\begin{proof}
Let $z_1>z_2>\cdots>z_n$ be the real roots of $P$ and~$z_1'>z_2'>\cdots>z_m'$ be the real roots of $Q$.
Observe that
\begin{gather*}
z_n=\min(Z_P)>\min(Z_Q)=z_m',\qquad
z_1=\max(Z_P)<\max(Z_Q)=z_1'.
\end{gather*}
Therefore, since the real roots of $P$ and~$Q$ interlace, we have
\begin{gather}
\label{rootsinterlace}
z_1'>z_1>z_2'>z_2>\cdots>z_n'>z_n>z_{n+1}'=z_m',
\end{gather}
In particular $m=n+1$.

Suppose $P(0)=0$.
Then there is an unique $1\leq k\leq n$ such that $z_k=0$.
So equation~\eqref{rootsinterlace} implies
\begin{gather*}
z_1'>z_1>z_2'>z_2>\!\cdots\!>z_{k-1}'\!>z_{k-1}\!>z_k'>z_k=0>z_{k+1}'\!>z_{k+1}\!>\!\cdots\!>z_n'>z_n>z_{n+1}'.
\end{gather*}
Therefore
\begin{gather*}
\left|Z_Q\cap (-\infty,0)\right|=n+1-(k+1)+1=n-k+1=\left|Z_P\cap (-\infty,0)\right|+1,\\
\left|Z_Q\cap (0,\infty)\right|=k=\left|Z_P\cap (0,\infty)\right|+1.
\end{gather*}
The case $Q(0)=0$ is proven similarly.

Suppose $P(0)>0$ and~$Q(0)>0$.
Since $P$ has a~positive leading coef\/f\/icient and~is not constant, we have
\begin{gather*}
\lim_{x\rightarrow\infty}P(x)=\infty.
\end{gather*}
Therefore, since $z_1$ is the largest real root of $P$, $P(x)>0$ for $x>z_1$.
Since $P$ has only simple roots, $P$ crosses $0$ positively at $z_1$, so $P(x)<0$ for $z_2<x<z_1$.
Again since $P$ has only simple roots, $P$ crosses $0$ negatively at $z_2$, so $P(x)>0$ for $z_3<x<z_2$.
Inductively we see that when $1\leq i<n$ is even, $P(x)>0$ for $z_{i+1}<x<z_i$, and~when $1\leq i<n$ is odd, $P(x)<0$ for $z_{i+1}<x<z_i$.
Furthermore $P(x)>0$ for $x<z_n$ if $n$ is even and~$P(x)<0$ for $x<z_n$ if $n$ is odd.

Similarly we have, for $1\leq i<n+1$ even, $Q(x)>0$ for $z_{i+1}'<x<z_i'$, and~for $1\leq i<n+1$ odd, $Q(x)<0$ for $z_{i+1}'<x<z_i'$.
Furthermore $Q(x)<0$ for $x<z_{n+1}'$, if $n$ is even and~$Q(x)>0$ for $x<z_{n+1}'$, if $n$ is odd.

There are three cases to consider: $z_1>0>z_n$, $z_1<0$ and~$z_n>0$.

We f\/irst assume $z_1>0>z_n$.
Then there is an unique $1\leq k\leq n$ such that $z_k>0>z_{k+1}$.
Since $z_k>0>z_{k+1}$ and~$P(0)>0$, we conclude that $k$ is even.
By equation~\eqref{rootsinterlace},
\begin{gather*}
z_k'>z_k>0>z_{k+1}>z_{k+2}'.
\end{gather*}
Since $k$ is even, $Q(x)>0$ for $z_{k+1}'<x<z_k'$ and~$Q(x)<0$ for $z_{k+2}'<x<z_{k+1}'$.
But $z_{k+2}'<0<z_k'$ and~$Q(0)>0$, hence $z_{k+1}'<0<z_k'$.
Therefore
\begin{gather*}
\left|Z_Q\cap (-\infty,0)\right|=n+1-(k+1)+1=\left|Z_P\cap (-\infty,0)\right|+1,\\
\left|Z_Q\cap (0,\infty)\right|=k=\left|Z_P\cap (0,\infty)\right|.
\end{gather*}
Let us assume $z_1<0$, then $P$ has no positive real roots.
Observe $Q(x)<0$ for $z_2'<x<z_1'$.
Suppose $z_1'>0$, then $z_2'>0$ since $Q(0)>0$.
Hence by equation~\eqref{rootsinterlace}, $z_1'>z_1>z_2'>0$, so $z_1>0$ and~we have a~contradiction.
So $z_1'<0$, hence all the real roots of $Q$ are negative and~we have
\begin{gather*}
\left|Z_Q\cap (-\infty,0)\right|=m=n+1=\left|Z_P\cap (-\infty,0)\right|+1,\\
\left|Z_Q\cap (0,\infty)\right|=0=\left|Z_P\cap (0,\infty)\right|.
\end{gather*}
Finally let us assume $z_n>0$, then $P$ has no negative real roots.
By equation~\eqref{rootsinterlace}, $z_n'>0$.
Since $P(0)>0$, $P(x)>0$ for $x<z_n$, therefore $n$ must be even.
Hence $Q(x)>0$ for $z_{n+1}'<x<z_n'$ and~$Q(x)<0$ for $x<z_{n+1}'$.
Since $z_n'>0$ and~$Q(0)>0$, this implies $z_{n+1}'<0<z_n'$.
Therefore
\begin{gather*}
\left|Z_Q\cap (-\infty,0)\right|=1=\left|Z_P\cap (-\infty,0)\right|+1,\qquad
\left|Z_Q\cap (0,\infty)\right|=n=\left|Z_P\cap (0,\infty)\right|.
\end{gather*}
This ends our discussion of the case $P(0)>0$ and~$Q(0)>0$.
The remaining cases are proven similarly.
\end{proof}

Taneda~\cite{taneda} proved that for $n\in\mathbb{N}$:
\begin{itemize}\itemsep=0pt
\item if $n\equiv1\pmod{3}$, then $\frac{Q_n}{z}\in\mathbb{Z}[z^3]$;
\item if $n\not\equiv1\pmod{3}$, then $Q_n\in\mathbb{Z}[z^3]$.
\end{itemize}
Hence $Q_n(0)=0$ if $n\equiv1\pmod{3}$.
By Theorem~\ref{simpleroots}, for every $n\geq1$, $Q_{n-1}$ and~$Q_n$ do not have a~common root.
Therefore $Q_n(0)=0$ if and~only if $n\equiv1\pmod{3}$.

Let us denote the coef\/f\/icient of the lowest degree term in $Q_n$ by $x_n$.
That is, we def\/ine $x_n:=Q_n(0)$ if $n\not\equiv1\pmod{3}$, and~$x_n:=Q_n'(0)$ if $n\equiv1\pmod{3}$.
In~\cite{roffelsen} we derived the following recursion for the $x_n$:
\begin{gather*}
x_0=1,\qquad x_1=1
\end{gather*}
and
\begin{gather}\label{EQ}
x_{n+1}x_{n-1}=
\begin{cases}(2n+1)x_n^2&\text{if $n\equiv0\pmod{3}$,}\\
4x_n^2 &\text{if $n\equiv1\pmod{3}$,}\\
-(2n+1)x_n^2 &\text{if $n\equiv2\pmod{3}$.}
\end{cases}
\end{gather}
We remark that the above recursion can be used to determine the $x_n$ explicitly, a~direct formula for~$x_n$ is given by Kaneko and~Ochiai~\cite{kaneko}.

\begin{lemma}
\label{lemsign}
For every $n\in\mathbb{N}$,
\begin{gather*}
\operatorname{sgn}(Q_n(0))=
\begin{cases}-1&\text{if $n\equiv3,5,6,8\pmod{12}$},\\
\hphantom{-}0 & \text{if $n\equiv1,4,7,10\pmod{12}$},\\
\hphantom{-}1 & \text{if $n\equiv0,2,9,11\pmod{12}$},
\end{cases}
\end{gather*}
where $\text{sgn}$ denotes the sign function on $\mathbb{R}$.
\end{lemma}
\begin{proof}
By induction using recursion~\eqref{EQ}, we have
\begin{gather*}
\operatorname{sgn}(x_n)=
\begin{cases}-1&\text{if $n\equiv3,5,7,6,8,10\pmod{12}$},\\
\hphantom{-}1 & \text{if $n\equiv0,1,2,4,9,11\pmod{12}$}.
\end{cases}
\end{gather*}
The lemma follows from this and~the fact that $Q_n(0)=0$ if and~only if $n\equiv1\pmod{3}$.
\end{proof}
We apply Lemma~\ref{lemnumbersroots} to the Yablonskii--Vorob'ev polynomials to prove Theorem~\ref{thmpositivenegative}.

\begin{proof}[Proof of Theorem~\ref{thmpositivenegative}.]
Let $n\geq2$, then by Proposition~\ref{proplimitbehaviour}, Theorem~\ref{simpleroots} and~Theorem~\ref{thminterlace}, $P:=Q_{n-1}$ and~$Q:=Q_{n+1}$ are monic polynomials with only simple roots such that the real roots interlace.
Furthermore by Theorem \ref{thmnumberrealroots},
both $P$ and~$Q$ have a~real root and
\begin{gather*}
\min(Z_P)>\min(Z_Q),\qquad
\max(Z_P)<\max(Z_Q).
\end{gather*}
So we can apply Lemma~\ref{lemnumbersroots} together with Lemma~\ref{lemsign} and~obtain:
\begin{gather*}
\left|Z_{n+1}\cap (-\infty,0)\right|=\left|Z_{n-1}\cap (-\infty,0)\right|+
\begin{cases} 0 & \text{if $n\equiv0,3\pmod{6}$,}\\
1 & \text{if $n\equiv1,2,4,5\pmod{6}$,}
\end{cases}\\
\left|Z_{n+1}\cap (0,\infty)\right|=\left|Z_{n-1}\cap (0,\infty)\right|+
\begin{cases} 0 & \text{if $n\equiv0,1\pmod{3}$,}\\
1 & \text{if $n\equiv2\pmod{3}$.}
\end{cases}
\end{gather*}
Observe that $Z_0=\varnothing$, $Z_1=\left\{0\right\}$ and~$Z_2=\left\{-\sqrt[3]{4}\right\}$.
The theorem is obtained by applying the above recursive formulas inductively.
\end{proof}

Let us discuss an example.
By Theorem~\ref{thmYV}, the unique rational solution of $P_{\rm II}(\alpha)$ for the parametervalue $\alpha:=21$ is given by
\begin{gather*}
w_{21}=\frac{Q_{20}'}{Q_{20}}-\frac{Q_{21}'}{Q_{21}}.
\end{gather*}
By Theorem~\ref{simpleroots}, $Q_{20}$ and~$Q_{21}$ do not have common roots and~the roots of $Q_{20}$ and~$Q_{21}$ are simple.
Hence the poles of $w_{21}$ are precisely the roots of $Q_{20}$ and~$Q_{21}$, the roots of $Q_{20}$ are poles of $w_{21}$ with residue $1$ and~the roots of $Q_{21}$ are poles of $w_{21}$ with residue $-1$.

By Theorem~\ref{thmnumberrealroots}, $Q_{20}$ has $10$ real roots and~by Theorem~\ref{thmpositivenegative}, $7$ of them are negative and~$3$ of them are positive.
Similarly $Q_{21}$ has $11$ real roots, $7$ of them are negative and~$4$ of them are positive.

Therefore $w_{21}$ has $21$ real poles, $10$ with residue~$1$ and~$11$ with residue~$-1$.
More precisely~$w_{21}$ has $7$ positive real poles, $3$ with residue~$1$ and~$4$ with residue~$-1$ and~$w_{21}$ has~$14$ negative real poles, $7$ with residue~$1$ and~$7$ with residue~$-1$.

\pdfbookmark[1]{References}{ref}

\LastPageEnding

\end{document}